\pdfoutput=1
\documentclass[AMA,STIX1COL]{WileyNJD-v2}

\usepackage[utf8]{inputenc}
\usepackage{microtype}
\usepackage{csquotes}

\usepackage{enumitem}
\setlist[enumerate,1]{label={(\roman*)}}

\usepackage{mathtools}
\usepackage{dsfont}

\newcommand*{\inputindex}{\omega}
\newcommand*{\consistencyspace}{\mathcal{X}}
\DeclarePairedDelimiter{\norm}{\lVert}{\rVert}
\DeclareMathOperator{\rk}{rk}
\DeclareMathOperator{\esssup}{ess\,sup}
\DeclareMathOperator{\im}{im}

\usepackage[nameinlink]{cleveref}
\crefname{eexample}{example}{examples}

\articletype{Research article}

\begin{document}

\title{Optimal control of differential-algebraic equations from an ordinary differential equation perspective}

\author[1]{Achim Ilchmann}
\author[1]{Leslie Leben}
\author[1]{Jonas Witschel*}
\author[1]{Karl Worthmann}

\authormark{ILCHMANN \textsc{et al.~}}

\address[1]{\orgdiv{Institut für Mathematik}, \orgname{Technische Universität Ilmenau}, \orgaddress{\state{98693 Ilmenau}, \country{Germany}}}

\corres{*Jonas Witschel, Technische Universität Ilmenau, Institut für Mathematik, Postfach 10 05 65, 98684 Ilmenau, Germany. \email{jonas.wischel@tu-ilmenau.de}}

\fundinginfohead{\textbf{Funding Information}}
\fundinginfoheadtext{}
\fundingInfo{The authors are indebted to the German Research Foundation (DFG) (grants IL 25/10-1 and WO2056/2-1) and the Studienstiftung des Deutschen Volkes for their support.}

\abstract[Summary]{
We study the Optimal Control Problem (OCP) for regular linear differential-algebraic systems (DAEs). To this end, we introduce the input index, which allows, on the one hand, to characterize the space of consistent initial values in terms of a Kalman-like matrix and, on the other hand, the necessary smoothness properties
of the control. The latter is essential to make the problem accessible from a numerical point of view. Moreover, we derive an augmented system as the key to analyze the OCP with tools well-known from optimal control of ordinary differential equations. The new concepts of the input index and the augmented system
  provide  easily checkable sufficient conditions which ensure that the stage costs are consistent with the differential-algebraic system.
}

\keywords{differential-algebraic equation, optimal control, consistent initial values, feedback control}

\jnlcitation{\cname{%
\author{Ilchmann A}, 
\author{Leben L}, 
\author{Witschel J}, and 
\author{Worthmann K}} (\cyear{2018}), 
\ctitle{Optimal control of differential-algebraic equations from an ordinary differential equation perspective}.}

\maketitle

\section{Introduction}
\noindent
We consider
time-invariant, single-input differential-algebraic systems (DAEs)
described, for
\([E,A,b]\in\mathbb{R}^{n\times n}\times\mathbb{R}^{n\times n}\times\mathbb{R}^n\) and $x^0\in\mathbb{R}^n$ ,
by
\begin{equation}
	\boxed{\tfrac{\mathrm{d}}{\mathrm{d} t} E x(t) = A x(t)+b u(t), \quad (E x)(0)=x^0,} \tag{DAE}\label{eq:dae}
\end{equation}
in \emph{quasi-Weierstraß form}, i.\,e.
\begin{equation*}
	E = \begin{bmatrix} I_{n_J} & 0 \\ 0 & N \end{bmatrix}, \qquad
	A = \begin{bmatrix} J & 0 \\ 0 & I_{n_N} \end{bmatrix}, \qquad
	b = \begin{pmatrix} b_J \\ b_N \end{pmatrix},
\end{equation*}
with \(J\in\mathbb{R}^{n_J\times n_J}\), \(b_J\in\mathbb{R}^{n_J}\), \(b_N\in\mathbb{R}^{n_N}\), and nilpotent \(N\in\mathbb{R}^{n_N \times n_N}\). %
It is well-known (Berger et al.\cite{BergIlch12a}, Th.~2.6) that every \emph{regular} DAE system, 
i.\,e. \(\det(s E-A)\neq0_{\mathbb{R}[s]}\), 
can be transformed into quasi-Weierstraß form, 
and that the dimension~\(n_N \in \mathbb{N}\) of the algebraic part is unique while the matrices \(J\) and~\(N\) are unique up to similarity. %
We assume $n_N \geq 1$, otherwise~\eqref{eq:dae} is an ordinary differential equation. For convenience, 
 \(x=(x_J,x_N)\) 
is partitioned according to the structure of the quasi-Weierstraß form.

The tuple~\((x,u)\in\mathcal{L}^1_\mathrm{loc}(\mathbb{R}_{\geq0},\mathbb{R}^n\times\mathbb{R})\) 
is called a \emph{solution} of~\eqref{eq:dae}  if, and only if,
\(E x\in\mathcal{W}^{1,1}_\mathrm{loc}(\mathbb{R}_{\geq0},\mathbb{R}^n)\) and \(u\) satisfy~\eqref{eq:dae} 
for almost all \(t\ge 0\).
An initial value $x^0\in\mathbb{R}^n$ is called \emph{consistent} if, and only if, a solution $(x,u)$
exists with $(E x)(0)=x^0$; the space of \emph{consistent initial values}  is denoted by~$\consistencyspace$.
\\

\noindent Essential for  our analysis will be the \emph{input index}~$\inputindex$ of~\eqref{eq:dae} defined by
\begin{equation}\label{eq:inputindex}
	\inputindex := \begin{cases} 
				\max \{i \in \mathbb{N} \mid N^i b_N \neq 0\} & \qquad \text{if $Nb_N \neq 0$}, \\
				0 &  \qquad \text{otherwise.}
	               \end{cases} 
\end{equation}
The input index is uniquely determined since~\(N\) and~\(b_N\) in the quasi-Weierstraß form
are unique up to similarity;
$\inputindex = 0$ holds if $b_N = 0$ or $Nb_N = 0$.
As  will be seen later, the input index is the number of times a function~$u$ must be differentiable almost everywhere to be a candidate for a solution $(x,u)$.
\\

The space~$\consistencyspace$ of consistent initial values will be characterized in 
terms of a Kalman like matrix and it will be shown that $\dim \consistencyspace = n_J + \inputindex$.
This sets us in a position to define an \emph{augmented system} of ordinary differential equations
which is equivalent to the~\eqref{eq:dae}.\\

Consider, for given weighting matrix \(S=S^\top\in\mathbb{R}^{(n+1)\times(n+1)}\) and consistent initial value~\(x^0\), the \emph{optimal control problem} over time horizon \([0,T)\), \(T\in(0,\infty]\):
\begin{equation}
	\boxed
	{
	\begin{aligned}
		\text{Minimize} \int_0^T \left(\begin{smallmatrix} x(t) \\ u(t) \end{smallmatrix}\right)^{\!\!\!\top} S 
		\left(\begin{smallmatrix} x(t) \\ u(t) \end{smallmatrix}\right) \mathrm{d} t 
		\quad\text{subject to}\quad (x,u) \text{ solves } \eqref{eq:dae}.
	\end{aligned}
	}
	\tag{OCP}\label{eq:daeocp}
\end{equation}
We determine the optimal value of the~\eqref{eq:daeocp} and an optimal feedback law by using the augmented system. 
This allows to use \emph{classical} ODE~results to solve~\eqref{eq:daeocp}.\\

The last aspect, namely to exploit \textit{classical} ODE~results in a DAE setting, is the key motivation of our approach since our long-term goal is the solution of (nonlinear) constrained optimal control problems on an infinite time horizon. %
However, constraints, e.\,g.\ mixed control-state constraints of the form $F x + g u \leq \mathds{1}$ with matrix $F \in \mathbb{R}^{p \times n}$ and vector $g \in \mathbb{R}^{p}$, in combination with an infinite optimization horizon render the problem, in general, computationally intractable.\cite{FalF14} %
Hence, we want to apply Model Predictive Control (MPC) to approximately solve this problem. %
In MPC, a sequence of \textit{constrained} optimal control problems on a finite time horizon~$T \in (0,\infty)$ is solved. %
While the theory is already well developed for ODE~systems, a rigorous stability analysis for the DAE~case is still missing. In addition, existing approaches do not preserve the particular structure of the underlying DAE.\cite{DiehBock02,bock2007constrained}

Essentially, there are two approaches in order to guarantee asymptotic stability of the origin w.r.t.\ the MPC closed loop for ODE-constrained systems. %
The first is based on stabilizing terminal constraints and costs.\cite{RawlingsMayneDiehl2017} Here, using the representation~\eqref{eq:optimalfeedback} of the optimal control allows to write the constraint as $(F + g\hat{k}_\alpha G^\dagger)x \leq \mathds{1}$ and to show existence of a weakly control invariant terminal set\cite{blanchini1999survey,grune2004asymptotic} $\{x^0 \in \mathcal{X} \mid V_\infty(x^0) \leq \rho\}$, $\rho > 0$.
Since the terminal set is a sub-level set of the value function~$V_\infty$, $V_\infty$ can be used to construct a suitable terminal cost. %
In conclusion, if Assumptions~\ref{ass:shatpossemidef} --- \ref{ass:shatnondegenerate} hold, recursive feasibility and asymptotic stability can be ensured by combining our results\cite{Ilch2018} with well-established techniques\cite{Mayne2000}. %
A similar approach\cite{Ilch2018b} --~again based on the results of this paper~-- works for~MPC without stabilizing terminal constraints and costs.\cite{GrunPann10,GrunPann11}
Furthermore, we are going to consider nonlinear differential-algebraic systems, whose linearization at the origin is stabilizable in a suitable sense.\cite{BergIlch12a} %
Again, \emph{standard} techniques for nonlinear ODEs should be applicable to exploit the presented results for a rigorous stability framework as it was done for ODEs.\cite{ChenAllgower1998,BocciaGrueneWorthmann2014}

The optimal control problem for DAEs has been studied by various authors. Basically, there exist three different approaches.\cite{Back06,KunkMehr08,ReisRend15} The first, which is also related to the work\cite{kurina1994control,kurina2004linear,LamoMarz13}, uses projectors. The second, which is further explicated in Reis and Voigt\cite{ReisVoig18}, uses the KYP inequality and Lur'e equations to characterize basic properties like feasibility and regularity of the~\eqref{eq:daeocp} with a zero terminal constraint at $T = \infty$ by \textit{storage functions} in the sense of Willems and Trentelmann\cite{willems2002synthesis,trentelman2002synthesis}. The third approach uses adjoint-based techniques and is applicable to nonlinear systems.\cite{biegler2012control} However, in none of them the optimal control~$u$ is explicitly specified, which is a basic requirement for the numerical treatment of optimal control problems governed by DAEs and, thus, for the usage of the \enquote{\emph{first discretize, then optimize}}-approach.\cite{gerdts2012optimal} But foremost, the main motivation to solve the optimal control problem via the quasi-Weierstraß form is that we consider this approach as very accessible from an ODE point of view. Hence, we can use the developed techniques to transfer well-known methods for ODE~systems to DAE-constrained systems as indicated above w.r.t.\ MPC. We plan to use this methodology as a blueprint to derive similar results based on the more sophisticated optimal control techniques \cite{Back06,KunkMehr08,ReisRend15} in a second step.

Overall, we propose a relatively simple approach to rewrite~-- under the condition that the DAE is already in 
quasi-Weierstraß form~-- \eqref{eq:daeocp} as an OCP constrained by \textit{ordinary} differential equations. The key ingredients are the \textit{input index} and the \textit{augmented system};
both of which are new concepts.

The structure of the paper is as follows. 
In~\Cref{SectionResultsDAE}, we exploit the quasi-Weierstraß form and show how much \enquote{freedom} in input function 
there is
and how the DAE is equivalent to an ODE.
In~\Cref{SectionResultsOCP},
the optimal control problem for~\eqref{eq:dae} is investigated. 
The optimal value function of the DAE is shown to coincide with an optimal value function of an
associated ODE optimal control problem.
Finally, in   \Cref{Sec:OCP-Feedback} it is proved how the optimal control can be expressed as a state feedback
for the DAE.

The Appendix consists of two parts. In Part~\ref{sec:odeoptimalcontrol}, classical results on optimal control of ODEs necessary
for our results are recalled. In Part~\ref{sec:proofs}, proofs of the results are given.
\\

\noindent\textbf{Notation}: The $(n \times n)$-identity matrix is denoted by \(I_n\). 
A matrix \(S\in\mathbb{R}^{n\times n}\) is called \emph{positive semi-definite}, denoted by \(S \succeq 0\), 
if, and only if, $x^\top S x \geq 0$ holds for all $x \in \mathbb{R}^n$;
\(N\in\mathbb{R}^{n\times n}\)
is called \emph{nilpotent} if, and only if, 
there exists $i \in \mathbb{N}$ such that $N^i=0$. 
The sets~\(\mathcal{L}^1_\mathrm{loc}(\mathbb{R}_{\geq0},\mathbb{R}^n)\) and \(\mathcal{W}^{k,1}_\mathrm{loc}(\mathbb{R}_{\geq0},\mathbb{R}^n)\) are the space of locally Lebesgue-integrable functions and the Sobolev space of functions whose derivatives up to the order~\(k\in\mathbb{N}_0\) exist almost everywhere and are in~\(\mathcal{L}^1_\mathrm{loc}(\mathbb{R}_{\geq0},\mathbb{R}^n)\), resp. The set $\mathcal{AC}(\mathbb{R}_{\geq0},\mathbb{R}^n)$ denotes the space of functions which are absolutely continuous   on each compact interval $K\subseteq\mathbb{R}_{\geq0}$.

\section{Differential-algebraic systems}\label{SectionResultsDAE}

\noindent As a convenient technical tool we define, for~\eqref{eq:dae} and the input index~$\inputindex$
as in~\eqref{eq:inputindex}, the \emph{Kalman like matrix}
\begin{equation}
	K := N \begin{bmatrix} b_N & \dots & N^{\inputindex-1} b_N \end{bmatrix} = \begin{bmatrix} N b_N & \dots & N^{\inputindex} b_N \end{bmatrix} \in \mathbb{R}^{n_N \times \inputindex} \label{eq:truncatedkalman}
\end{equation}
with the convention that $K$ is void if $\inputindex = 0$.

Any solution of~\eqref{eq:dae} can be characterized in terms of the input index
and the Kalman like matrix as follows.

\begin{proposition}\label{thm:differentiability}
Consider~\eqref{eq:dae}. Then $(x,u)$ is a solution to~\eqref{eq:dae} if, and only if,
\begin{subequations}
\begin{align}
\label{eq:solxn-a}
	u(\cdot) \in & \ \mathcal{W}^{\inputindex,1}_\mathrm{loc}(\mathbb{R}_{\geq 0},\mathbb{R}),  \\
\label{eq:solxn-b}
	x_J(t) = & \ \mathrm{e}^{Jt} x_J(0) + \int_0^t \mathrm{e}^{J(t-s)} b_J u(s) \mathrm{d} s  \qquad  \text{ for all } t \geq 0,\\
\nonumber
	x_N(t) = & \ - \sum_{i=0}^{\inputindex} N^i b_N u^{(i)}(t) 
	= - \left[ b_N,\, K\right] (u(t), \ldots, u^{(\inputindex)}(t))^\top\\
\label{eq:solxn-c}	
	& \hspace{44mm} \quad\text{for almost all } t\geq0.
\end{align}
\end{subequations} 
\end{proposition}
\noindent
The proof is given in~\Cref{sec:proofs}.

Note that the sum in~\eqref{eq:solxn-c} goes up to~$\inputindex$ which may be smaller than~$n_N$
as stated in Berger and Reis\cite{BergReis13a}, Remark~3.4.

\Cref{thm:differentiability} implies
in particular that not every initial value $x^0\in \mathbb{R}^n$ is
consistent.
In the following proposition, which is proved in~\Cref{sec:proofs}, some properties and a characterization of the space of consistent initial data
\(\consistencyspace\subseteq\mathbb{R}^n\)
are presented.

\begin{proposition}\label{thm:consistentiv}
For~\eqref{eq:dae}, we have:
\begin{enumerate}
\item\label{item:kalmanrank}
The matrix \(K\) defined by~\eqref{eq:truncatedkalman} has full rank, i.\,e. $\rk K = \inputindex  \leq n_N$.
\item\label{item:consistencyspace}
$\consistencyspace=\im F \quad\text{ with }\quad  
  F:=\begin{bmatrix} I_{n_J} & 0 \\ 0 & -K \end{bmatrix} \in \mathbb{R}^{n \times (n_J+\inputindex)}$.
\item\label{item:leftinverse}
\(F\) has a left inverse \(F^\dagger\in\mathbb{R}^{(n_J+\inputindex) \times n}\).
\item\label{item:dimconsistencyspace}
$\dim \consistencyspace = n_J+\inputindex$.
\end{enumerate}
\end{proposition}

As an immediate consequence of \Cref{thm:differentiability} and \Cref{thm:consistentiv}
we collect the following corollary.

\begin{corollary}\label{Cor:consistentiv}
If  $(x,u)$ is a solution to~\eqref{eq:dae},
then
\[
F \begin{pmatrix} x_J(0) \\ u(0) \\ \vdots \\ u^{(\inputindex-1)}(0) \end{pmatrix} = (Ex)(0)= x^0 \in \consistencyspace\,.
\]
\end{corollary}

\ \\

For ODEs, i.e.\ $n_N=0$ in~\eqref{eq:dae}, it follows from~\eqref{eq:solxn-b}
that for any initial value 
$x^0\in \mathbb{R}^n$ and any input function
$u\in \mathcal{L}^1_\mathrm{loc}(\mathbb{R}_{\geq0},\mathbb{R})$
there exists a unique solution $(x,u)$ of~\eqref{eq:dae}.
In other words, the mapping
\[
\mathbb{R}^n \times \mathcal{L}^1_\mathrm{loc}(\mathbb{R}_{\geq0},\mathbb{R}) \to 
 \mathcal{W}^{1,1}_\mathrm{loc}(\mathbb{R}_{\geq0},\mathbb{R}^n),
 \
 t\mapsto x(t)= \mathrm{e}^{Jt} x_J(0) + \int_0^t \mathrm{e}^{J(t-s)} b_J u(s) \, \mathrm{d} s
\] 
is well-defined.
Thus $(x^0,u)$
can be \enquote{freely} chosen and determines a unique~$x$
so that
$(x,u)$ solves the initial value problem~\eqref{eq:dae}.
In case of~\eqref{eq:dae} with $n_N\ge 1$, formula~\eqref{eq:solxn-c}
shows that~$x$ and~$u$ must satisfy an algebraic equation. Hence, %
$x^0$ and $u$
cannot be chosen freely. 
To be precise, we define the set
\[
\mathcal{F} := 
\left\{ (x^0,u) \in \consistencyspace \times \mathcal{W}^{\inputindex,1}_\mathrm{loc}(\mathbb{R}_{\geq0},\mathbb{R})
\,\middle|\,
F \begin{pmatrix} x_J(0) \\ u(0) \\ \vdots \\ u^{(\inputindex-1)}(0) \end{pmatrix} =  x^0
\right\}
\]
and the mapping
\begin{equation}
\begin{array}{rccl}
\varphi : & \mathcal{F} & \to & \mathcal{L}^1_\mathrm{loc}(\mathbb{R}_{\geq0},\mathbb{R}^n)\\[1ex]
          &  (x^0,u)    & \mapsto & 
  x=   \begin{pmatrix} x_J\\ x_N  \end{pmatrix}  =
    \left( \begin{smallmatrix} \mathrm{e}^{J \cdot} x_J(0) + \int_0^\cdot \mathrm{e}^{J(\cdot -s)} b_J u(s) \mathrm{d} s\\[1ex]
    - \left[ b_N,\, K\right] (u(\cdot), \ldots, u^{(\inputindex)}(\cdot))^\top \end{smallmatrix}\right).
\end{array}
\end{equation}
Then, for any $(x^0,u)\in \mathcal{F}$, the tuple $(\varphi(x^0,u),u)$ solves the initial value problem~\eqref{eq:dae}.

However, there is some freedom in~$u$ left as the following remark shows.

\begin{remark}\label{Prop:freedom-u}
Consider~\eqref{eq:dae} and define
\begin{equation*}
\mathcal{U}^0 := 
\bigl\{ u\in  \mathcal{W}^{\inputindex,1}_\mathrm{loc}(\mathbb{R}_{\geq0},\mathbb{R})
\bigm|
u(0)= \dots = u^{(\inputindex-1)}(0)  =   0
\bigr\}.
\end{equation*}
Then for any 
$(x^0,u)  \in \consistencyspace \times \mathcal{W}^{\inputindex,1}_\mathrm{loc}(\mathbb{R}_{\geq0},\mathbb{R})$
we have 
\[
 (x^0,u)  \in \mathcal{F}
 \quad \Longleftrightarrow
 \quad 
 \forall \widetilde u\in \mathcal{U}^0 : \
  (x^0,u+ \widetilde u)  \in \mathcal{F}.
\]
\ \\[-5ex]
\ ~\hspace*{3cm}\hfill $\diamond$
\end{remark}

\Cref{Prop:freedom-u} shows that whenever $(x,u)$ is a solution to~\eqref{eq:dae},
then for any $\widetilde u\in \mathcal{U}^0$ there exists a solution for $(x(0),u+\widetilde u)$.

\begin{remark}
If $\inputindex = 0$, then \(K\) is void, \(F\) is given by 
\(\begin{bmatrix} I_{n_J} & 0_{n_J\times n_N} \end{bmatrix}^\top\) and \Cref{thm:consistentiv} states that $x^0 \in \consistencyspace$ holds if, and only if, $x^0_i = 0$ for all $i \in \{n_J+1,\ldots,n\}$. %
Hence, \Cref{thm:differentiability} shows that \eqref{eq:dae} reduces to an initial value problem for the ordinary differential equation
\begin{equation}\label{eq:augmentedsystem:one}
	\dot{x}_J(t)=J x_J(t)+b_J u(t)\text{,}
\end{equation}
with initial value \(x_J(0)=F^\dagger x^0\) and \enquote{output} \(x_N(t)=-b_N u(t)\).
\hfill $\diamond$
\end{remark}

For $\inputindex \geq 1$, the above findings allow to rewrite~\eqref{eq:dae} equivalently as an ordinary differential equation
where 
lower-order derivatives are introduced as new system states.

\begin{proposition}\label{thm:augmentedsystem}
Consider the initial value problem~\eqref{eq:dae} with $Nb_N \neq 0$,
let the input index~$\inputindex$ be given by~\eqref{eq:inputindex}
and define 
\begin{equation*}
	\hat{x} = \begin{pmatrix} x_J^\top, & u, & \dots, & u^{(\inputindex-1)} \end{pmatrix}^{\!\!\top},
	\qquad
	\hat{u}=u^{(\inputindex)},
		\quad\text{and}\quad
	\hat{n} = n_J + \inputindex.
\end{equation*}
Then \((x,u)\) with \((E x)(0)=x^0\in\consistencyspace\) is a solution to~\eqref{eq:dae} if, and only if, \((\hat{x},\hat{u},x_N)\) with \(\hat{x}(0)=F^\dagger x^0\)  is a solution to the \emph{augmented system}
\begin{subequations}\label{eq:augmentedsystem}
\begin{align}
	\dot{\hat{x}}(t) &= \underbrace{\begin{bmatrix} J & b_J & 0 & \dots & 0 \\ 0 & 0 & 1 & \ddots & \vdots \\ \vdots & \ddots & \ddots & \ddots & 0 \\ \vdots & & \ddots & 0 & 1 \\ 0 & \dots & \dots & 0 & 0 \end{bmatrix}}_{=:\hat{A}\, \in\, \mathbb{R}^{\hat{n}\times\hat{n}}} \hat{x}(t)+\underbrace{\begin{pmatrix} 0 \\ \vdots \\ \vdots \\ 0 \\ 1 \end{pmatrix}}_{=:\hat{b}\, \in\, \mathbb{R}^{\hat{n}}} \hat{u}(t) \label{eq:xo} \\
	x_N(t) &= -\begin{bmatrix} 0_{n_N\times n_J} & b_N & K \end{bmatrix} \begin{pmatrix} \hat{x}(t) \\ \hat{u}(t) \end{pmatrix}.
\end{align}
\end{subequations}

\begin{proof}
Using \Cref{thm:differentiability} and reading~\eqref{eq:dae} and system~\eqref{eq:augmentedsystem} row-wise, 
yields that \((x,u)\) solves~\eqref{eq:dae} if, and only if, \((\hat{x},\hat{u},x_N)\) solves~\eqref{eq:augmentedsystem}.
The correspondence of the initial values follows from
\begin{align*}
	F \hat{x}(0)
	&=  \begin{bmatrix} I_{n_J} & 0 \\ 0 & -\begin{bmatrix} N b_N & \dots & N^{\inputindex} b_N \end{bmatrix} \end{bmatrix} 
	\begin{pmatrix} x_J(0) \\ u(0) \\ \vdots \\ u^{\inputindex-1}(0) \end{pmatrix} \\
	&= \begin{pmatrix} x_J (0) \\ -\sum\limits_{i=1}^{\inputindex} N^i b_N u^{(i-1)}(0) \end{pmatrix}
	= \begin{pmatrix} x_J(0) \\ N x_N(0) \end{pmatrix}
	= (E x)(0).
\end{align*}
\end{proof}
\end{proposition}

The advantage of the augmented system~\eqref{eq:augmentedsystem}
 is that its input~\(\hat{u}\), unlike the input of~\eqref{eq:dae}, can be chosen freely from \(\mathcal{L}^1_\mathrm{loc}(\mathbb{R}_{\geq0},\mathbb{R})\). %
Furthermore, it inherits stabilizability properties of~\eqref{eq:dae} which will  
 be useful for the optimal control problem considered in the next section.
\begin{definition} \label{Def:beh_stab}
\eqref{eq:dae} is called \emph{behaviourally stabilizable} if, and only if,
\begin{align*}
	\begin{lgathered}
	\forall (x,u) \text{ solution of \eqref{eq:dae}}\,\exists \widetilde{u}\in\mathcal{U}^0: \\
	\lim_{t\to\infty} \esssup_{[t,\infty)} \norm{(\varphi((E x)(0),u+\widetilde{u}),u+\widetilde{u})}=0.
	\end{lgathered}
\end{align*}
\end{definition}
In \Cref{Prop:stab} behavioural stabilization will be characterized algebraically
by the Hautus criterion.
This result then shows by~Berger and Reis\cite{BergReis13a}, Corollary~4.3
that \Cref{Def:beh_stab} is equivalent to
the standard definition:
\begin{align*}
	\begin{lgathered}
	\forall (x,u) \text{ solution of~\eqref{eq:dae}}\,\exists (\tilde{x},\tilde{u})\text{ solution of~\eqref{eq:dae}}: \\
	(x,u)|_{(-\infty,0)} \overset{\mathrm{ae}}{=} (\tilde{x},\tilde{u})|_{(-\infty,0)} \land {} \lim_{t\to\infty} \esssup_{[t,\infty)} \lVert(\tilde{x},\tilde{u})\rVert=0.
	\end{lgathered}
\end{align*}
Note that in the former the time axis is $\mathbb{R}_{\geq0}$, whereas in Berger and Reis\cite{BergReis13a} it
is $\mathbb{R}$; however, for time-invariant systems this does not make a difference.

\begin{proposition}\label{Prop:stab}
\eqref{eq:dae} is behaviourally stabilizable if, and only if,
\begin{equation*}
	\forall \lambda\in\overline{\mathbb{C}}_+: \mathrm{rk}_\mathbb{C}  [\lambda E-A,b] = \rk_{\mathbb{R}(s)} [s E-A,b].
\end{equation*}

\begin{proof}
As \(\det(s E-A)\neq0_{\mathbb{R}[s]}\), it holds that \(\rk_{\mathbb{R}(s)} [s E-A,b]=n\). Furthermore as~\eqref{eq:dae} is in quasi-Weierstraß form and \(\forall \lambda\in\mathbb{C}: \rk (\lambda N-I)=n_N\), we have that
\begin{align}
	&\mathrel{\hphantom{\Longleftrightarrow}} \forall \lambda\in\overline{\mathbb{C}}_+: \mathrm{rk}_\mathbb{C}  [\lambda E-A,b] = \rk_{\mathbb{R}(s)} [s E-A,b]=n \notag \\
	&\Longleftrightarrow \forall \lambda\in\overline{\mathbb{C}}_+: \mathrm{rk}_\mathbb{C}  [\lambda I_{n_J}-J,b_J] = n_J. \label{eq:odestabilizable}
\end{align}
By~\eqref{eq:solxn-b} and~\Cref{thm:consistentiv}\ref{item:consistencyspace}, it follows that~\eqref{eq:odestabilizable} is a necessary condition for the stabilizability of the DAE. 
For   sufficiency, suppose that~\eqref{eq:odestabilizable} holds. For a solution~\((x,u)\) of~\eqref{eq:dae}, choose any~\(\widetilde{u}\in\mathcal{U}^0\) so that~\(\bar{u}:=u+\widetilde{u}\) stabilizes the ODE part of~\eqref{eq:dae} given by~\eqref{eq:augmentedsystem:one}. By~\Cref{thm:differentiability}, it follows that \((\varphi((E x)(0),\bar{u}),\bar{u})\) converges to zero almost everywhere.
\end{proof}
\end{proposition}

\begin{lemma}\label{thm:stabilizability}
Consider~\eqref{eq:dae} with \(\inputindex>0\). Then~\eqref{eq:dae} is behaviourally stabilizable if, and only if, the ODE system~\((\hat{A},\hat{b})\) given by the augmented system~\eqref{eq:augmentedsystem} is stabilizable.

\begin{proof}
By \Cref{Prop:stab} and the Hautus criterion we have that
\begin{align*}
	&\mathrel{\hphantom{\Longleftrightarrow}} \text{\eqref{eq:dae} is behaviourally stabilizable} \\
	&\Longleftrightarrow  \forall \lambda\in\overline{\mathbb{C}}_+: \rk_\mathbb{C} [\lambda E-A,b] = \rk_{\mathbb{R}(s)} [s E-A,b] \\
	&\Longleftrightarrow \forall \lambda\in\overline{\mathbb{C}}_+: \rk_\mathbb{C}[\lambda I_{n_J}-J,b_J] = n_J \\
	&\Longleftrightarrow \forall \lambda\in\overline{\mathbb{C}}_+:
	\mathrm{rk}_\mathbb{C}[\lambda I-\hat{A},\hat{b}] = \rk_\mathbb{C} \begin{bmatrix} \lambda I_{n_J}-J & b_J & 0 & \dots & 0 \\ 0 & 0 & 1 & \ddots & \vdots \\ \vdots & \ddots & \ddots & \ddots & 0 \\ 0 & \dots & 0 & 0 & 1 \end{bmatrix} = \hat{n} \\
	&\Longleftrightarrow (\hat{A},\hat{b}) \text{ is stabilizable}. \qedhere
\end{align*}
\end{proof}
\end{lemma}

\begin{eexample}\label{example:exampleaugmented}
To illustrate the previous results, consider the DAE
\begin{equation}
	\frac{\mathrm{d}}{\mathrm{d} t}
	\begin{bmatrix} 0 & 1 & 0 \\ 0 & 0 & 1 \\ 0 & 0 & 0 \end{bmatrix}
	x(t)
	=
	x(t)
	+
	\begin{pmatrix} 0 \\ 1 \\ 0 \end{pmatrix}
	u(t),
	\quad
	\left(\begin{bmatrix} 0 & 1 & 0 \\ 0 & 0 & 1 \\ 0 & 0 & 0 \end{bmatrix} x\right)(0)=x^0.
	\label{eq:exampledae}
\end{equation}
It has an input index~\(\inputindex=2\) and the consistent initial values are, due to~\Cref{thm:consistentiv},
given by
\begin{equation*}
	\consistencyspace = \im F = \im (-K) = \im \begin{pmatrix} -1 \\ 0 \\ 0 \end{pmatrix} .
\end{equation*}
\Cref{thm:augmentedsystem} yields the augmented system
\begin{align*}
	\dot{\hat{x}}(t) &= \hat{u}(t), \\
	x(t) &= -\begin{bmatrix} 0 & 1 \\ 1 & 0 \\ 0 & 0 \end{bmatrix} \begin{pmatrix} \hat{x}(t) \\ \hat{u}(t) \end{pmatrix} = \begin{pmatrix} -\hat{u}(t) \\ -\hat{x}(t) \\ 0 \end{pmatrix},
\end{align*}
where \(\hat{x}=u\), \(\hat{u}=\dot{u}\), and with an initial value of
\begin{equation*}
	\hat{x}(0) = F^\dagger x^0 = (-1,0,0) \, x^0 = -x^0_1.
\end{equation*}
\ \\[-5ex]
\ ~\hspace*{3cm}\hfill $\diamond$
\end{eexample}

\section{Optimal control}\label{SectionResultsOCP}
\noindent
For a given weighting matrix \(S=S^\top\in\mathbb{R}^{(n+1)\times(n+1)}\) and time horizon   $T\in(0,\infty]$, 
the \emph{cost functional} $J_T $
assigns to each solution~\((x,u)\) of~\eqref{eq:dae} the value
\begin{align}\label{eq:costfunctional}
	J_T(x,u) = \int_0^T \begin{psmallmatrix} x(t) \\ u(t) \end{psmallmatrix}^{\!\!\!\top} S 
		\, \begin{psmallmatrix} x(t) \\ u(t) \end{psmallmatrix} \, \mathrm{d} t
		\in  \mathbb{R}\cup\{\pm\infty\} .
\end{align}
Our goal is to find the infimum of the cost functional for given consistent initial value \(x^0\in\consistencyspace\), denoted by the \emph{optimal value function}
\begin{equation}\label{eq:VT}
	\begin{aligned}
V_T : \consistencyspace  &\to  \mathbb{R}\cup\{\pm\infty\} \\
        x^0   &\mapsto \inf_{(x,u)} J(x,u) \quad\text{subject to $(x,u)$ solves~\eqref{eq:dae}}.
        \end{aligned}
\end{equation}
\ 
\\

We define the optimal value function associated to the augmented system~\eqref{eq:augmentedsystem}
in a similar fashion.
The above weighting matrix~$S$   is partitioned according to
 the structure of the quasi-Weierstra\ss\ form~\eqref{eq:dae} as
\begin{equation*}
	S =: \begin{bmatrix} Q_J & Q_{JN} & h_J \\ Q_{JN}^\top & Q_N & h_N \\ h_J^\top & h_N^\top & r\end{bmatrix},
\qquad
\text{where \(Q_J\in\mathbb{R}^{n_J\times n_J}\), \(Q_N\in\mathbb{R}^{n_N\times n_N}\), \(r\in\mathbb{R}\)},
\end{equation*}
and a new symmetric matrix is defined as
	\begin{equation}
		\hat{S} = \begin{bmatrix} 
			Q_J & h_J-Q_{JN} b_N & -Q_{JN} K \\ 
			h_J^\top-b_N^\top Q_{JN}^\top & r+b_N^\top Q_N b_N-2 h_N^\top b_N & (b_N^\top Q_N-h_n^\top) K \\
			-K^\top Q_{JN}^\top & K^\top (Q_N b_N-h_n) & K^\top Q_N K 
		\end{bmatrix}
		\in \mathbb{R}^{(\hat{n}+1)\times(\hat{n}+1)} \label{eq:shat}
	\end{equation}
where $\hat{n} = n_J + \inputindex$. If $\inputindex = 0$, then~$K$ is void and the matrix~$\hat{S}$ reduces to a $(2 \times 2)$-block matrix.
The \emph{cost functional}  $\hat J_T$
assigns to each solution~\((\hat x,\hat u)\) of~\eqref{eq:augmentedsystem} the value
\begin{align}\label{eq:costfunctionalaugmented}
	\hat J_T(\hat x,\hat u) = \int_0^T \begin{psmallmatrix} \hat x(t) \\ \hat u(t) \end{psmallmatrix}^{\!\!\!\top} \hat S 
		\, \begin{psmallmatrix} \hat x(t) \\ \hat u(t) \end{psmallmatrix} \, \mathrm{d} t
		\in  \mathbb{R}\cup\{\pm\infty\} .
\end{align}
Finally,  the \emph{optimal value function}
is
\begin{equation*}
	\begin{aligned}
\hat V_T : \consistencyspace  &\to  \mathbb{R}\cup\{\pm\infty\} \\
        x^0   &\mapsto \inf_{(\hat x, \hat u)} \hat J(\hat x, \hat u) \quad\text{subject to $(\hat x,\hat u)$ 
        solves~~\eqref{eq:augmentedsystem} and $ \hat{x}(0)=F^\dagger x^0$}.
        \end{aligned}
\end{equation*}

We will show in the next theorem that  the DAE optimal control problem~\eqref{eq:daeocp} is equivalent 
to the ODE optimal control problem
\begin{equation*}
	\boxed
	{
	\begin{aligned}
&\text{Minimize } \hat J_T(\hat x,\hat u)
	\quad\text{subject to $\hat{u} \in \mathcal{L}_{\textrm{loc}}^1 \text{, } \hat{x}(0)=F^\dagger x^0$} \\
	& \phantom{\text{Minimize }}  \text{and the ODE } \begin{cases}  
		\eqref{eq:xo} & \text{for $\inputindex \geq 1$,} \\
		\eqref{eq:augmentedsystem:one} & \text{for $\inputindex = 0$.}
	\end{cases} \nonumber
	\end{aligned}
	}
	\tag{ODE-OCP}\label{eq:equivalentocp}
\end{equation*}
Equivalent means that the optimal values are the same for a given consistent initial value~$x^0$.

\begin{theorem}\label{thm:optimalvalue}
The optimal values of the optimal control problems~\eqref{eq:daeocp} and~\eqref{eq:equivalentocp} coincide,
i.e.
	\begin{align}
		V_T(x^0) = \hat{V}_T(x^0) \qquad
		\text{for any consistent initial value~\(x^0 \in \consistencyspace\).}
		\nonumber
	\end{align}
\end{theorem}
\noindent
The theorem is proved in~\Cref{sec:proofs}.
 
\Cref{thm:optimalvalue} allows to apply well-known results of optimal control of ODEs
and, for example, to derive a representation of the optimal value function~$V_T$ in terms of a quadratic form and give a sufficient condition for feasibility of~\eqref{eq:daeocp}, i.\,e.~\(V_T\) is finite on $\consistencyspace$. 
To define the differential Riccati equation, we partition the symmetric weighting matrix~$\hat{S}$ given by~\eqref{eq:shat}  as
\begin{equation*}
	\hat{S} = \begin{bmatrix*}[l] \hat{Q} & \hat{h} \\ \hat{h}^\top & \hat{r} \end{bmatrix*} \quad\text{with}\quad \hat{Q} \in \mathbb{R}^{\hat{n} \times \hat{n}} \text{ and } \hat{r} \in \mathbb{R}
\end{equation*}
and the relevant   \emph{differential Riccati equation} becomes
\begin{equation}
	\begin{aligned}
	\dot{P}(t) = \hat{A}^\top P(t)+P(t) \hat{A}+\hat{Q}-(P(t) \hat{b}+\hat{h}) \hat{r}^{-1} (P(t) \hat{b}+\hat{h})^{\hspace*{-.75mm}\top} \hspace*{-.75mm}\text{, }\quad P(0) = 0,
	\end{aligned}
	\tag{DRE}\label{eq:riccatiode}
\end{equation}
where \(\hat{A}\) and \(\hat{b}\) are given by~\eqref{eq:xo} (or \(\hat{A}=J\), \(\hat{b}=b_J\) if \(\inputindex=0\)). 
The limit \(\lim_{t\to\infty} P(t)\) exists, cf.~\Cref{sec:odeoptimalcontrol}, and is denoted by~\(P(\infty)\).
\\

The following assumptions are formulated in terms of  \(\hat{S}\)  defined by~\eqref{eq:shat} and $\hat{A}$ given as in~\eqref{eq:xo} for $\inputindex \geq 1$ or $J$ for $\inputindex = 0$:
\hspace*{0.75cm}\begin{enumerate}[label={\textbf{(A\arabic*)}},ref={(A\arabic*)},leftmargin=*]
	\item\label{ass:shatpossemidef} \(\hat{S}\succeq0\),
	\item\label{ass:stabilizable} \eqref{eq:dae} is stabilizable, i.\,e. the pair \((J,b_J)\) is stabilizable,
	\item\label{ass:rhatpos} 
\(\hat{r} =
\begin{pmatrix} 0 & \dots & 0 & 1 \end{pmatrix} \hat{S} \begin{pmatrix} 0 &\dots & 0 & 1 \end{pmatrix}^\top  >0\),
	\item\label{ass:observable} 
\(\bigl( \hat{A}, \hat{Q} \bigr) =
 \bigl( \hat{A}, \begin{bmatrix} I_{\hat{n}} & 0_{\hat{n}\times 1} \end{bmatrix} \hat{S} \begin{bmatrix} I_{\hat{n}} & 0_{\hat{n}\times 1} \end{bmatrix}^\top \bigr)  \) \ is observable,
	\item\label{ass:shatnondegenerate} 
\(\rk \hat{S}
 = \rk \hat{Q} + 1
 = \rk \bigl(\begin{bmatrix} I_{\hat{n}} & 0_{\hat{n}\times1} \end{bmatrix} \hat{S} \begin{bmatrix} I_{\hat{n}} & 0_{\hat{n}\times1} \end{bmatrix}^\top\bigr)+1 \).
\end{enumerate}

\noindent
The first three assumptions imply that the value function~$V_T$ can be represented by a quadratic function, thus, exhibits finite values on~$\mathcal{X}$. If, in addition, also the forth and fifth assumption hold, the optimal value function is positive definite.

\begin{theorem}\label{thm:optimalvalueproperties}
Suppose \ref{ass:shatpossemidef}\,--\,\ref{ass:rhatpos} hold and let~$P(T)$ and~$V_T$
 be defined by~\eqref{eq:riccatiode} and~\eqref{eq:VT}, resp. Then
the following holds for	\(T\in(0,\infty]\) :
	\begin{enumerate}
		\item $\forall\, x^0\in\consistencyspace: 
V_T(x^0)=(x^0)^\top  (F^\dagger)^\top \, P(T)\, F^\dagger x^0$ and \(P(T)\succeq0\).
		\item If, in addition,  \ref{ass:observable} and~\ref{ass:shatnondegenerate} hold, then \(P(T)\succ0\), i.\,e. \\
		\(\forall x^0\in\consistencyspace\setminus\{0\}: V_T(x^0)>0\).
	\end{enumerate}

	\begin{proof}
		By \Cref{thm:stabilizability}, Assumption~\ref{ass:stabilizable} is equivalent to the stabilizability of the augmented system~\eqref{eq:xo}. %
		Therefore, the ODE~\eqref{eq:xo} (or~\eqref{eq:augmentedsystem:one} for \(\inputindex=0\)) fulfils the assumptions of~\Cref{thm:optimalvalueode}. %
		By~\Cref{thm:optimalvalue}, the ODE optimal control problem in~\Cref{thm:optimalvalue} is equivalent to~\eqref{eq:daeocp}, thus showing the assertion.
	\end{proof}
\end{theorem}

%
\section{The  optimal control as a feedback}\label{Sec:OCP-Feedback}
\noindent
It is well-known that in case of ODEs the optimal control can be expressed as a state feedback.
We will show that a similar result holds for~\eqref{eq:dae}. For \(\inputindex=0\), this is immediately clear from \Cref{thm:optimalcontrol} as \(x^*=\hat{x}^*\) holds. For \(\inputindex>0\), we need some further algebraic manipulations to obtain an optimal feedback. 
We stress that in case of DAEs the closed-loop system is not necessarily regular;
this is shown in the following example.

\begin{eexample}
Consider the scalar DAE \(0=x(t)+u(t)\), which has the form~\eqref{eq:dae}, where \(n_J=0\), \(N=0\), \(b_N=1\). Choosing the feedback \(u(t)=-x(t)\), we obtain the closed-loop system \(0=x(t)-x(t)=0\). Hency any function~\(x\in\mathcal{L}^1_\mathrm{loc}(\mathbb{R}_{\geq0},\mathbb{R})\) is a solution of the closed-loop system, which is not regular any more. This happens because the feedback only provides the superfluous  information \(u=-x\), which is already containted in the algebraic constrains of the DAE.
\hfill $\diamond$
\end{eexample}

The following ~\namecref{thm:singularfeedback}   characterizes when a feedback leads to a regular closed-loop system.

\begin{lemma}\label{thm:singularfeedback}
Consider~\eqref{eq:dae} and \(k=(k_J,k_N)\in\mathbb{R}^{1\times n_J}\times\mathbb{R}^{1\times n_N}\). Then
\begin{equation}
	\begin{aligned}
	& \det(s E-(A-b k))=0_{\mathbb{R}[s]} \\
	\Longleftrightarrow\ {} & k_J (s I_{n_J}-J)^{-1} b_J=0 \ \land\  k_N b_N = 1 \ \land\ k_N K=0_{1\times\inputindex}. \label{eq:feedbacksingular}
	\end{aligned}
\end{equation}
\end{lemma}
The proof is given  in~\Cref{sec:proofs}.

\ \\

\ \\

Consider the (OCP) for~\eqref{eq:dae}
with solution $(x^*,u^*)$ and consistent initial value~$x^0$.
Then~$x^*$ and~$u^*$ are called \textit{optimal state trajectory} and \emph{optimal control}
for time horizon $T\in(0,\infty]$ if, and only if,
$J_T(x^*,u^*)= V_T(x^0)$.
\\

We show that the optimal control of the~(OCP) equals 
the ($n_J+1$)-component of the augmented state~\(\hat{x}\),
where the latter solves an ODE.

\begin{proposition}\label{thm:optimalcontrol}
Consider the optimal control problem~\eqref{eq:daeocp} with optimization horizon $T \in (0,\infty]$ and
suppose that Assumptions~\ref{ass:shatpossemidef}\,--\,\ref{ass:rhatpos} hold. Then, the unique optimal control
is given by $u^*=\hat{x}^*_{n_J+1}(t)$, where $\hat{x}^*$ solves the ODE
\begin{equation}
	\dot{\hat{x}}^*(t) = \bigl[\hat{A}-\hat{b} \, \hat{r}^{-1} (\hat{b}^\top P(T-t)+\hat{h}^\top)\bigr] \hat{x}^*(t), \quad \hat{x}^*(0) = F^\dagger x^0, \label{eq:optimaltrajectoryaugmented}
\end{equation}
and \(P\) is the solution of the differential Riccati equation~\eqref{eq:riccatiode}. 
If $\inputindex = 0$, then $u^*(t) = -\hat{r}^{-1} (\hat{b}^\top P(T-t)+\hat{h}^\top) \hat{x}^*(t)$.

\begin{proof}
Let \(x^0\in\consistencyspace\) be arbitrary, and consider the ODE optimal control problem for the augmented system~\eqref{eq:xo} stated in~\Cref{thm:optimalvalue}. 
In view of \Cref{thm:optimalvalueode},
the optimal trajectory~\(\hat{x}^*\) and optimal control~\(\hat{u}^*\) for the augmented system are given by the solution of~\eqref{eq:optimaltrajectoryaugmented} and \(\hat{u}^*(t)=-\hat{r}^{-1} (\hat{b}^\top P(T-t)+\hat{h}^\top) \hat{x}^*(t)\), resp. Therefore~\Cref{thm:optimalvalue} yields
\begin{equation*}
	V_T(x^0) = \int_0^T \begin{psmallmatrix} \hat{x}^*(t) \\ \hat{u}^*(t) \end{psmallmatrix}^{\!\!\top} \hat{S}\ \begin{psmallmatrix} \hat{x}^*(t) \\ \hat{u}^*(t) \end{psmallmatrix}.
\end{equation*}
By~\Cref{thm:augmentedsystem} and~\Cref{thm:optimalvalue},
\begin{equation*}
	\begin{pmatrix} x_J^* \\ x_N^* \\ u^* \end{pmatrix}
	=
	\begin{bmatrix} I_{n_J} & 0 & 0 \\ 0 & -b_N & -K \\ 0 & 1 & 0 \end{bmatrix} \begin{pmatrix} \hat{x}^* \\ \hat{u}^* \end{pmatrix}
\end{equation*}
is a solution of~\eqref{eq:dae} that fulfils
\begin{equation*}
	J_T(x^*,u^*) = \int_0^T \begin{psmallmatrix} \hat{x}^*(t) \\ \hat{u}^*(t) \end{psmallmatrix}^{\!\!\top} \hat{S}\ \begin{psmallmatrix} \hat{x}^*(t) \\ \hat{u}^*(t) \end{psmallmatrix} = V_T(x^0).
\end{equation*}
Therefore \(u^*\) is an optimal control.

To prove the uniqueness, assume that we have a solution \((x',u')\) of~\eqref{eq:dae} with \(J_T(x',u')=J_T(x^*,u^*)\). 
\Cref{thm:augmentedsystem} and \Cref{thm:optimalvalue} yield that
\begin{equation*}
	\hat{x}' = \begin{pmatrix} {x'_J}^\top, & u', & \dots, & {u'}^{(\inputindex-1)} \end{pmatrix}^\top, \quad
	\hat{u}' = {u'}^{(\inputindex)}
\end{equation*}
is an optimal solution of the ODE optimal control problem for the augmented system~\eqref{eq:xo}. 
Thus uniqueness of the ODE optimal control (cf.~\Cref{thm:optimalvalueode}\ref{thm:optimalcontrolode})
yields \(\hat{u}'=\hat{u}^*\) and hence \(\hat{x}'=\hat{x}^*\). For \(\inputindex=0\), this immediately shows the uniqueness of~\(u^*\) as \(u'=\hat{u}'=\hat{u}^*=u^*\). For \(\inputindex>0\), it follows from~\Cref{thm:augmentedsystem} that
\begin{equation*}
	u'= \hat{x}'_{n_J+1} = \hat{x}^*_{n_J+1} = u^*,
\end{equation*}
which concludes the proof.
\end{proof}
\end{proposition}
In the infinite optimization horizon case \(T=\infty\), note that by convention~\(P(\infty-t)=P(\infty)\) for all \(t\geq0\), so that~\eqref{eq:optimaltrajectoryaugmented} is a time-invariant ODE.

\begin{eexample} \label{Ex:3}
Revisit~\Cref{example:exampleaugmented} and define for the DAE~\eqref{eq:exampledae}
the cost functional
\begin{equation*}
	J_T(x,u) = \int_0^T \norm{x(t)}^2+u(t)^2 \,\mathrm{d} t.
\end{equation*}
To solve the optimal control problem, rewrite the~(OCP) as an ODE optimal control problem for the augmented system as 
described  in~\Cref{thm:optimalvalue}: The cost functional for the augmented system is
\begin{equation*}
	\hat{J}_T(\hat{x},\hat{u}) = \int_0^T \begin{pmatrix} \hat{x}(t) \\ \hat{u}(t) \end{pmatrix}^\top \begin{pmatrix} 2 & 0 \\ 0 & 1 \end{pmatrix} \begin{pmatrix} \hat{x}(t) \\ \hat{u}(t) \end{pmatrix},
\end{equation*}
so the Riccati equation~\eqref{eq:riccatiode} takes the form
\[
	\dot{P}(t) = 2-P(t)^2, \quad P(0)=0,
\]
and hence
\[
\forall T\geq0: P(T)=\frac{\sqrt{2} \, (\mathrm{e}^{2\sqrt{2} T}-1)}{\mathrm{e}^{2\sqrt{2} T}+1}, \  P(\infty)=\lim_{t\to\infty} P(t)=\sqrt{2}.
\]
\Cref{thm:optimalvalue} yields the optimal value
\begin{equation*}
	V_T(x^0) = \hat{V}_T(x^0) = (x^0)^\top (F^\dagger)^\top P(T) F^\dagger x^0 = P(T) \cdot (x^0_1)^2.
\end{equation*}
In view of~\Cref{thm:optimalcontrol}, the optimal control for the DAE~\eqref{eq:exampledae} satisfies
the ODE initial value problem
\begin{equation}
	\dot{u}^*(t) = -P(T-t) u^*(t), \quad u^*(0)=-x^0_1. \label{eq:exampleoptimalcontrolaugmented}
\end{equation}
In case of an infinite time horizon \(T=\infty\), this leads to the optimal trajectory and control given for almost all \(t\geq0\) by
\begin{equation}
	x^*(t) = \begin{psmallmatrix} -\sqrt{2} \mathrm{e}^{-\sqrt{2} t} x_1^0 \\ \mathrm{e}^{-\sqrt{2} t} x_1^0 \\ 0 \end{psmallmatrix}, \quad
	u^*(t) = -\mathrm{e}^{-\sqrt{2} t} x_1^0. \label{eq:exampleoptimaltrajectory}
\end{equation}
\hfill $\diamond$
\end{eexample}

We are now in a position to state the main result of this section,
that is how the optimal control can be expressed as a feedback which yields
a regular closed-loop system. 

\begin{theorem}\label{thm:optimalfeedback}
Consider~\eqref{eq:daeocp} with optimization horizon $T = \infty$ and $Nb_N \neq 0$ and suppose that Assumptions~\ref{ass:shatpossemidef}\,--\,\ref{ass:rhatpos} hold. Define
\begin{equation*}
	\hat{k}_\alpha 
	:= \alpha \begin{pmatrix} \hat{r}^{-1} (\hat{b}^\top P(\infty) +\hat{h}^\top) & 1 \end{pmatrix}+\begin{pmatrix} 0_{1\times n_J} & 1 & 0_{1\times\inputindex} \end{pmatrix} \in \mathbb{R}^{1\times(\hat{n}+1)},
	\quad
\alpha \in \mathbb{R},
\end{equation*}
where \(P(\infty) = \lim_{t \to \infty} P(t) \in \mathbb{R}^{\hat{n} \times \hat{n}}\) and~\(P\) 
is the solution of the Riccati equation~\eqref{eq:riccatiode}. Furthermore let \(G^\dagger\in\mathbb{R}^{(\hat{n}+1)\times n}\) be a left inverse of
\begin{equation*}
	G = \begin{bmatrix} I_{n_J} & 0 & 0 \\ 0 & -b_N & -K \end{bmatrix}.
\end{equation*}
Then for any \(\alpha\neq0\), the closed-loop system
\begin{equation}
	\tfrac{\mathrm{d}}{\mathrm{d} t} E x(t) = (A+b \hat{k}_\alpha G^\dagger) \, x(t), \quad (E x)(0)=x^0 \label{eq:optimalclosedloop}
\end{equation}
obtained by the feedback
\begin{equation}
	u(t)=\hat{k}_\alpha \, G^\dagger x(t) \label{eq:optimalfeedback}
\end{equation}
is regular and its solution is an optimal state trajectory for~\eqref{eq:daeocp}.
\end{theorem}
The proof is carried out in~\Cref{sec:proofs}.

\begin{remark}
While the optimal \emph{control} of~\eqref{eq:dae} is unique, the optimal \emph{feedback} is generally not unique: any \(\alpha\neq0\) and any left inverse of \(G\) can be chosen for  the feedback, giving rise to a possible multitude of optimal feedbacks. Roughly speaking, this is due to the fact that the states~\(u\), \dots, \(u^{(\inputindex-1)}\) of the augmented system which are necessary for the optimal control can be derived from~\eqref{eq:dae} in multiple equivalent ways. Nevertheless all feedbacks will lead to the same optimal control input.
\hfill $\diamond$
\end{remark}

\begin{eexample} \label{Ex:4}
We continue with \Cref{Ex:3}.
For \(T=\infty\), the optimal control can be expressed as a state feedback: 
Using the family of left inverses
\begin{equation*}
	G^\dagger = \begin{bmatrix} 0 & -1 & 0 \\  -1 & 0 & \beta \end{bmatrix}, \quad \beta\in\mathbb{R},
\end{equation*}
we arrive by~\Cref{thm:optimalfeedback} at the family of distinct optimal feedbacks
\begin{equation*}
	u = \begin{pmatrix} -\alpha & -\sqrt{2}\alpha-1 & \beta \end{pmatrix} x, \quad \alpha\in\mathbb{R}\setminus\{0\}, \ \beta\in\mathbb{R},
\end{equation*}
leading to the closed-loop system
\begin{equation*}
	\frac{\mathrm{d}}{\mathrm{d} t}
	\begin{bmatrix} 0 & 1 & 0 \\ 0 & 0 & 1 \\ 0 & 0 & 0 \end{bmatrix}
	x(t)
	=
	\begin{bmatrix} 1 & 0 & 0 \\ -\alpha & -\sqrt{2} \alpha & \beta \\ 0 & 0 & 1 \end{bmatrix} x(t), \quad
	\left(\begin{bmatrix} 0 & 1 & 0 \\ 0 & 0 & 1 \\ 0 & 0 & 0 \end{bmatrix} x\right)(0)=x^0.
\end{equation*}
The solution of the closed-loop system fulfils for almost all \(t\geq0\)
\begin{equation*}
	x(t) = \begin{psmallmatrix} -\sqrt{2} \mathrm{e}^{-\sqrt{2} t} x_1^0 \\ \mathrm{e}^{-\sqrt{2} t} x_1^0 \\ 0 \end{psmallmatrix},
\end{equation*}
coinciding with the optimal trajectory stated in \eqref{eq:exampleoptimaltrajectory}.
\hfill $\diamond$
\end{eexample}

\section{Conclusions and Outlook}\label{SectionConclusions}

Most of the presented results should be directly generalizable to multi-input DAEs by defining the input index as a vector, whose entries correspond to the columns of the matrix~$B$. Furthermore, we conjecture that similar results can also be established for non-regular DAEs.

The presented results on optimal control of regular DAEs may be worth knowing in its own right. However, our guided research interest stems from Model Predictive Control (MPC) with finite time horizon~$T \in (0,\infty)$ in order to approximately solve~\eqref{eq:daeocp} with mixed control-state constraints
\begin{equation}\label{eq:controlstateconstraint}
	F x + g u \leq \mathds{1}\text{, }\quad (F,g) \in \mathbb{R}^{p \times n}  \times \mathbb{R}^{p}
\end{equation}
on an infinite time horizon. As outlined in the introduction, this hope is justified, 
see~Ilchmann et al.\cite{Ilch2018,Ilch2018b} for details. In particular, our approach helps to identify the control and also to characterize initial feasibility in terms of the consistency space. Moreover, also a suitable combination of terminal region and terminal costs was deduced in Ilchmann et al.\cite{Ilch2018} such that the origin is asymptotically stable w.r.t.\ the MPC closed loop. %
In conclusion, the proposed approach allows to directly transfer \textit{standard} techniques used in MPC to systems governed by DAEs.

The choice of the quasi-Weierstra\ss{} form in~\eqref{eq:dae} is mainly motivated by the fact that it makes DAE optimal control very accessible from an ODE point of view. %
In the future we want to investigate whether similar results are achievable based on the Feedback Equivalence Form (FEF)\cite{ReisRend15} or representation concepts introduced in~Embree and Blake.\cite{embree2017pseudospectra} %
The FEF is numerically more robust in comparison to the quasi-Weierstra\ss{} form, but may also allow for the application of \textit{classical} ODE~results since the lower right block of the transformed matrix~$A$ is invertible. Hence, we think that the presented approach is the starting point to design MPC schemes based on more sophisticated optimal control techniques.\cite{Back06,KunkMehr08,ReisRend15}

\appendix

\section{Optimal control of ordinary differential equations}\label{sec:odeoptimalcontrol}

\noindent
Consider for 
\([\bar{A},\bar{b}]\in\mathbb{R}^{\bar{n}\times\bar{n}}\times\mathbb{R}^{\bar{n}}\)
and  \(x^0\in\mathbb{R}^{\bar{n}}\)
the initial value problem of ordinary differential equations 
\begin{equation}
	\dot{\bar{x}}(t) = \bar{A} \bar{x}(t) + \bar{b} \bar{u}(t), \quad \bar{x}(0)=\bar{x}^0 \label{eq:ode}
\end{equation}
with associated cost functional
\begin{equation*}
	\bar{J}_T(\bar{x},\bar{u}) = \int_0^T \begin{psmallmatrix} \bar{x}(t) \\ \bar{u}(t) \end{psmallmatrix}^{\!\!\top} \bar{S} \ \begin{psmallmatrix} \bar{x}(t) \\ \bar{u}(t) \end{psmallmatrix} \mathrm{d} t,
\end{equation*}
where 
\(\bar{S}=\bar{S}^\top\in\mathbb{R}^{(\bar{n}+1)\times(\bar{n}+1)}\), and optimal control problem
\begin{equation*}
	\begin{aligned}
	\bar{V}_T : \mathbb{R}^{\bar{n}}  &\to  \mathbb{R}\cup\{\pm\infty\} \\
        x^0   &\mapsto \inf_{(\bar{x}, \bar{u})} \bar{J}(\bar{x}, \bar{u}) \quad\text{subject to $(\bar{x},\bar{u})$ 
        solves~~\eqref{eq:ode} and } \bar{x}(0)=\bar{x}^0.
        \end{aligned}
\end{equation*}

To guarantee the feasibility of the optimal control problem, i.\,e. to ensure that \(\bar{V}_T(\bar{x}^0)\) is finite for all \(T\in(0,\infty]\)
and all  \(\bar{x}^0\in\mathbb{R}^n\), we introduce some standard assumptions.\cite{LancRodm95} These are already formalized in~\ref{ass:shatpossemidef}\,--\,\ref{ass:shatnondegenerate} for the \emph{DAE} optimal control problem~\eqref{eq:daeocp} that are identical to the assumptions for the ODE optimal control problem, except that in~\ref{ass:stabilizable}, the stabilizability of~\eqref{eq:dae} needs to be replaced by the stabilizability of the ODE~\eqref{eq:ode}.

\begin{proposition}\label{thm:riccatiode}
Let the ODE~\eqref{eq:ode} be stabilizable and Assumptions~\ref{ass:shatpossemidef} and~\ref{ass:rhatpos} be fulfilled. Then the associated differential Riccati equation~\eqref{eq:riccatiode} satisfies
\begin{enumerate}
\item there exists a solution~\(P(\cdot): [0,\omega)\to\mathbb{R}^{n\times n}\) for some maximal \(\omega\in(0,\infty]\),
\item the life span of this solution can be extended on the whole half axis, i.\,e. \(\omega=\infty\),
\item the solution is unique, symmetric and positive semidefinite, i.\,e. \(P(t)=P(t)^\top\succeq0\) for all \(t>0\),
\item the solution is monotonically non-decreasing in time, i.\,e.
\begin{equation*}
	\forall \delta>0\,\forall t\geq0: P(t+\delta)\succeq P(t),
\end{equation*}
\item\label{item:are} \(\exists P(\infty)=P(\infty)^\top\in\mathbb{R}^{n\times n}: \lim_{t\to\infty} P(t)=P(\infty)\succeq0\).
\end{enumerate}
If additionally Assumptions~\ref{ass:observable} and \ref{ass:shatnondegenerate} hold, then
\begin{enumerate}[resume]
\item \(\forall t>0: P(t)\succ0\), \(P(\infty)\succ0\).
\end{enumerate}

\begin{proof}
According to Lancaster and Rodman\cite{LancRodm95}, Th.~16.4.3, the differential Riccati equation~\eqref{eq:riccatiode} has a unique solution on~\([0,\infty)\) and this solution is monotonically non-decreasing. Using Lancaster and Rodman\cite{LancRodm95}, Prop.~16.2.5 and~Th.~16.4.4, we obtain that \(P(\infty):=\lim_{t\to\infty} P(t)\) exists.

If Assumptions~\ref{ass:observable} and \ref{ass:shatnondegenerate} hold, \(P(\infty)\) is positive definite according to Lancaster and Rodman\cite{LancRodm95}, Th.~16.3.3. It remains to be shown that \(P(t)\succ0\) for all \(t>0\). This follows analogously to~Lancaster and Rodman\cite{LancRodm95}, Prop.~16.2.8.
\end{proof}
\end{proposition}

\begin{remark}
According to Lancaster and Rodman\cite{LancRodm95}, Th.~16.4.4 and Lem.~10.11, the limit \(P(\infty)\) in~\Cref{thm:riccatiode}\ref{item:are} is the minimal solution of the \emph{algebraic Riccati equation}
\begin{equation}
	0 = \bar{A}^\top P+P \bar{A}+\bar{Q}-(P \bar{b}+\bar{h}) \bar{r}^{-1} (P \bar{b}+\bar{h})^{\hspace*{-.75mm}\top}, \tag{ARE}\label{eq:are}
\end{equation}
i.\,e. for all \(P\in\mathbb{R}^{\bar{n}\times\bar{n}}\) that fulfil~\eqref{eq:are}, it holds that \(P-P(\infty)\succeq0\).
\end{remark}

We are now in a position to explain the consequences of the differential Riccati equation~\eqref{eq:riccatiode} when applied to the optimal control problem.

\begin{proposition}\label{thm:optimalvalueode}
Let the ODE~\eqref{eq:ode} be stabilizable and suppose the associated optimal control problem fulfils Assumptions~\ref{ass:shatpossemidef} and~\ref{ass:rhatpos}. Then we have, for \(T\in(0,\infty]\):
\begin{enumerate}
\item There exists a unique \(P_T=P_T^\top\succeq0\) such that
\begin{equation*}
	\forall \bar{x}^0\in\mathbb{R}^{\bar{n}}: \bar{V}_T(\bar{x}^0)=(\bar{x}^0)^\top P_T \bar{x}^0
\end{equation*}
and this \(P_T\) is given by \(P_T=P(T)\)  in \Cref{thm:riccatiode}.
\item\label{thm:optimalcontrolode} For every \(\bar{x}^0\in\mathbb{R}^{\bar{n}}\), the unique optimal control \(\bar{u}^*\in\mathcal{L}^1_\mathrm{loc}([0,T),\mathbb{R})\) such that \(V_T(\bar{x}^0)=\bar{J}_T(\bar{x}^0,\bar{u}^*)\) is given by
\begin{align*}
	\dot{\bar{x}}^*(t) &= \bigl[\bar{A}-\bar{b} \, \bar{r}^{-1} (\bar{b}^\top P(T-t)+\bar{h}^\top)\bigr] \bar{x}^*(t), \quad \bar{x}^*(0) = \bar{x}^0 \\
	\bar{u}^*(t) &= -\bar{r}^{-1} (\bar{b}^\top P(T-t)+\bar{h}^\top) \bar{x}^*(t),
\end{align*}
where \(P(\cdot)\) is the solution of~\eqref{eq:riccatiode} and by convention \(P(\infty-t)=P(\infty)\) for all \(t\geq0\).
\end{enumerate}
If additionally Assumptions~\ref{ass:observable} and \ref{ass:shatnondegenerate} hold, then
\begin{enumerate}[resume]
\item \(\forall T\in(0,\infty]\,\forall \bar{x}^0\in\mathbb{R}^{\bar{n}}: \bar{V}_T(\bar{x}^0)>0\).
\end{enumerate}

\begin{proof}
This follows according to Lancaster and Rodman{LancRodm95}, Lem.~16.4.2 and Th.~16.3.3 for \(T\in(0,\infty)\) and \(T=\infty\), resp.
\end{proof}
\end{proposition}

\section{Proofs}\label{sec:proofs}

\begin{proof}[\textbf{Proof of~\Cref{thm:differentiability}}]
\eqref{eq:solxn-b} follows immediately from the variation of constants formula. For~\eqref{eq:solxn-a} and~\eqref{eq:solxn-c}, consider the algebraic part of~\eqref{eq:dae}, i.\,e.
\begin{equation}
	\tfrac{\mathrm{d}}{\mathrm{d} t} N x_N(t) = x_N(t)+b_N u(t). \label{eq:ndae}
\end{equation}
In passing, we note that \(\varphi\in\mathcal{W}^{1,1}_\mathrm{loc}(\mathbb{R}_{\geq0},\mathbb{R}^n)\) if, and only if, \(\varphi_1,\dots,\varphi_n\in\mathcal{W}^{1,1}_\mathrm{loc}(\mathbb{R}_{\geq0},\mathbb{R})\) and therefore \(M \varphi\in\mathcal{W}^{1,1}_\mathrm{loc}(\mathbb{R}_{\geq0},\mathbb{R}^l)\) for any \(M\in\mathbb{R}^{l\times n}\).

Let \((\begin{psmallmatrix} x_J \\ x_N \end{psmallmatrix},u\bigr)\) be a solution of~\eqref{eq:dae}, then \((x_N,u)\) is a solution of~\eqref{eq:ndae}. We proceed in several steps.

\emph{Step 1}: We show that
\begin{equation}
	\forall i>\inputindex: N^i x_N=0. \label{eq:nixn0}
\end{equation}
For \(i\geq\min\{i\in\mathbb{N} \mid N^i=0\}>\inputindex\), we have \(N^i=0\) and therefore~\eqref{eq:nixn0} follows.

Assume that~\eqref{eq:nixn0} holds for \(i>\inputindex+1\). Then we have
\begin{equation*}
	N^{i-1} x_N = N^{i-1} (x_N+b_N u) \underset{\text{\eqref{eq:ndae}}}{\overset{\mathrm{ae}}{=}} N^{i-1} \tfrac{\mathrm{d}}{\mathrm{d} t} N x_N = \tfrac{\mathrm{d}}{\mathrm{d} t} N^i x_N \overset{\text{\eqref{eq:nixn0}}}{=} 0.
\end{equation*}
Therefore~\eqref{eq:nixn0} is shown for~\(i-1\).

\emph{Step 2}: For \(\inputindex=0\), we have \(x_N=\frac{\mathrm{d}}{\mathrm{d} t}(N x_N)-b_N u=-b_N u\) and~\eqref{eq:solxn-c} is shown. For \(\inputindex>0\), we show the following statement by induction:
\begin{align}
	\forall i\in\{0,\dots,\inputindex-1\}:{} & u\in \mathcal{W}^{i+1,1}_\mathrm{loc}(\mathbb{R}_{\geq0},\mathbb{R}) \text{ and} \label{eq:uw12} \\
	& N^{\inputindex} b_N u^{(i)}
	\overset{\mathrm{ae}}{=} -N^{\inputindex-i} x_N
	- \sum_{k=0}^{i-1} N^{\inputindex-i+k} b_N u^{(k)}. \label{eq:diffNbu}
\end{align}
For \(i=0\), we have
\begin{gather*}
	0 \overset{\text{\eqref{eq:nixn0}}}{=} \tfrac{\mathrm{d}}{\mathrm{d} t} N^{\inputindex+1} x_N = N^{\inputindex} \tfrac{\mathrm{d}}{\mathrm{d} t} N x_N \overset{\mathrm{ae}}{\underset{\eqref{eq:ndae}}{=}} N^{\inputindex} x_N+N^{\inputindex} b_N u \\
	\Longrightarrow\ N^{\inputindex} b_N u \overset{\mathrm{ae}}{=} -N^{\inputindex} x_N \in \mathcal{W}^{1,1}_\mathrm{loc}(\mathbb{R}_{\geq0},\mathbb{R}^n)
\end{gather*}
and so~\eqref{eq:diffNbu} follows. Furthermore, as \(N^{\inputindex} b_N\neq0\) and \(u\) is scalar, we have~\eqref{eq:uw12}.

Assume that~\eqref{eq:uw12} and~\eqref{eq:diffNbu} hold for \(i\in\{0,1,\dots,\inputindex-2\}\). Then we have that
\begin{align}
	\tfrac{\mathrm{d}}{\mathrm{d} t} N^{\inputindex} b_N u^{(i)}
	&\underset{\mathclap{\eqref{eq:diffNbu}}}{\overset{\mathrm{ae}}{=}} \frac{\mathrm{d}}{\mathrm{d} t} \Bigl(-N^{\inputindex-i} x_N - \sum_{k=0}^{i-1} N^{\inputindex-i+k} b_N u^{(k)}\Bigr) \notag \\
	&\overset{\mathrm{ae}}{\underset{\mathclap{\eqref{eq:ndae}}}{=}} -N^{\inputindex-(i+1)} (x_N+b_N u) - \sum_{k=0}^{i-1} N^{\inputindex-i+k} b_N u^{(k+1)} \notag \\
	&= -\underbrace{N^{\inputindex-(i+1)} x_N}_{\in\mathcal{W}^{1,1}_\mathrm{loc}(\mathbb{R}_{\geq0},\mathbb{R}^n)} - \sum_{k=0}^{(i+1)-1} \underbrace{N^{\inputindex-(i+1)+k} b_N u^{(k)}}_{\overset{\eqref{eq:uw12}}\in\mathcal{W}^{1,1}_\mathrm{loc}(\mathbb{R}_{\geq0},\mathbb{R}^n)}. \label{eq:diffNbuinduction}
\end{align}
This shows~\eqref{eq:diffNbu} for \(i+1\). 

Furthermore, as the right-hand side of~\eqref{eq:diffNbuinduction} is in \(\mathcal{W}^{1,1}_\mathrm{loc}(\mathbb{R}_{\geq0},\mathbb{R}^n)\), \(N^{\inputindex} b_N\neq0\) and \(u\) is scalar, it follows that \(u^{(i+1)}\in\mathcal{W}^{1,1}_\mathrm{loc}(\mathbb{R}_{\geq0},\mathbb{R})\). Therefore~\eqref{eq:uw12} holds for \(i+1\).

\emph{Step 3}: For \(\inputindex>0\), reconsider~\eqref{eq:diffNbu} for \(i=\inputindex-1\): by~\eqref{eq:uw12}, we know that \(u\in\mathcal{W}^{i+1,1}(\mathbb{R}_{\geq0},\mathbb{R})\), so~\eqref{eq:diffNbu} is differentiable. Carrying out the differentiation in the same way as in~\eqref{eq:diffNbuinduction}, we arrive at
\begin{equation}
	N^{\inputindex} b_N u^{(\inputindex)} \overset{\mathrm{ae}}{=} -x_N - \sum_{k=0}^{\inputindex-1} N^k b_N u^{(k)}. \label{eq:diffNbuind-1}
\end{equation}
Rearranging~\eqref{eq:diffNbuind-1} for \(x_N\) immediately gives~\eqref{eq:solxn-c}.
\end{proof}

\begin{proof}[\textbf{Proof of~\Cref{thm:consistentiv}}]
For~\ref{item:kalmanrank}, we show by induction for \(0\leq i\leq\inputindex\) that
\begin{equation}\label{eq:truncatedkalmanrank}
	\rk \begin{bmatrix} N^i b_N, & N^{i+1} b_N & \dots & N^{\inputindex} b_N \end{bmatrix} = \inputindex-i+1.
\end{equation}
For \(i=\inputindex\), \eqref{eq:truncatedkalmanrank} follows from \(N^{\inputindex} b_N\neq0\).

Assume that~\eqref{eq:truncatedkalmanrank} holds for an \(i\in\{1,\dots,\inputindex\}\). Consider the linear combination
\begin{equation}\label{eq:linearcombination}
	\alpha_{i-1} N^{i-1} b_N+\dots+\alpha_{\inputindex} N^{\inputindex} b_N = 0 
\end{equation}
for \(\alpha_{i-1},\dots,\alpha_{\inputindex}\in\mathbb{R}\). Multiplying~\eqref{eq:linearcombination} from the left by~\(N\), we arrive at
\begin{align*}
	0
	&= \alpha_{i-1} N^i b_N+\dots+\alpha_{\inputindex-1} N^{\inputindex} b_N+\alpha_{\inputindex} N^{\inputindex+1} b_N \\
	&= \alpha_{i-1} N^i b_N+\dots+\alpha_{\inputindex-1} N^{\inputindex} b_N.
\end{align*}
By \eqref{eq:truncatedkalmanrank}, we see that
\begin{equation*}
	\alpha_{i-1} = \dots = \alpha_{\inputindex-1} = 0.
\end{equation*}
Therefore by \eqref{eq:linearcombination} and \(N^{\inputindex} b_N\neq0\), we obtain \(\alpha_{\inputindex}=0\) and hence \eqref{eq:truncatedkalmanrank} for \(i-1\). Now \eqref{eq:truncatedkalmanrank} for \(i=1\) shows the assertion.

For~\ref{item:consistencyspace}, note that by \Cref{thm:differentiability}, it follows that every solution~\((x,u)\) of~\eqref{eq:dae} fulfils
\begin{equation*}
	(E x)(0) = \begin{bmatrix} x_J(0) \\ N x_N(0) \end{bmatrix} = \begin{bmatrix} x_J(0) \\ -\sum_{i=1}^{\inputindex} N^i b_N u^{(i-1)}(0) \end{bmatrix} = \begin{bmatrix} I_{n_J} & 0 \\ 0 & -K \end{bmatrix} \begin{pmatrix} x_J(0) \\ u(0) \\ \vdots \\ u^{(\inputindex-1)}(0) \end{pmatrix},
\end{equation*}
where \(x_J(0)\), \(u(0)\), \dots, \(u^{\inputindex-1}(0)\) can be chosen freely.

Now~\ref{item:leftinverse} and~\ref{item:dimconsistencyspace} immediately follow from~\ref{item:kalmanrank} and~\ref{item:consistencyspace}.
\end{proof}

\begin{proof}[\textbf{Proof of~\Cref{thm:optimalvalue}}]
As a first step, we show that for any solution~\((x,u)\) of~\eqref{eq:dae}, the cost functional~\eqref{eq:costfunctional} is given by
\begin{equation}
	J_T(x,u) = \hat{J}_T(\hat{x},\hat{u}), \label{eq:costfunctionalequality}
\end{equation}
where \((\hat{x},\hat{u})\) is the solution of the augmented system~\eqref{eq:xo} (or of~\eqref{eq:augmentedsystem:one} for \(\inputindex=0\)) with \(\hat{x}(0)=F^\dagger x^0\): by~\Cref{thm:augmentedsystem} (or by~\(x_N=-b_N u=-b_N \hat{u}\) for \(\inputindex=0\)), it follows that
\begin{equation*}
	\begin{pmatrix} x_J \\ x_N \\ u \end{pmatrix}
	\overset{\mathrm{ae}}{=}
	\begin{bmatrix} I_{n_J} & 0 & 0 \\ 0 & -b_N & -K \\ 0 & 1 & 0 \end{bmatrix} \begin{pmatrix} \hat{x} \\ \hat{u} \end{pmatrix}.
\end{equation*}
Substituting this in~\eqref{eq:costfunctional} gives~\eqref{eq:costfunctionalequality} as
\begin{equation*}
	\begin{bmatrix} I_{n_J} & 0 & 0 \\ 0 & -b_N & -K \\ 0 & 1 & 0 \end{bmatrix}^\top
	\begin{bmatrix} Q_J & Q_{JN} & h_J \\ Q_{JN}^\top & Q_N & h_N \\ h_J^\top & h_N^\top & r \end{bmatrix}
	\begin{bmatrix} I_{n_J} & 0 & 0 \\ 0 & -b_N & -K \\ 0 & 1 & 0 \end{bmatrix}
	=
	\hat{S}.
\end{equation*}

Now let \(x^0\in\consistencyspace\) and \(\varepsilon>0\) be arbitrary. By \(\hat{V}_T: \consistencyspace\to\mathbb{R}\cup\{\pm\infty\}\) we denote the optimal value function of the ODE optimal control problem
\begin{equation*}
		\hat{V}_T(x^0) = \inf_{\hat{u} \in \mathcal{L}_{\mathrm{loc}}^1} \hat{J}_T(\hat{x},\hat{u})  \quad\text{subject to}\quad \text{\eqref{eq:xo} and } \hat{x}(0)=F^\dagger x^0
\end{equation*}
(or~\eqref{eq:augmentedsystem:one} instead of~\eqref{eq:xo} for \(\inputindex=0\)). Choose \(\hat{u}\in\mathcal{L}^1_\mathrm{loc}(\mathbb{R},\mathbb{R})\) such that \(\hat{J}_T(\hat{x},\hat{u})\leq \hat{V}_T(x^0)+\varepsilon\), where~\(\hat{x}\) is a solution of the ordinary differential equation~\eqref{eq:xo} (or~\eqref{eq:augmentedsystem:one} for \(\inputindex=0\)) with \(\hat{x}(0)=F^\dagger x^0\). Define
\begin{equation*}
	\begin{pmatrix} x_J \\ x_N \\ u \end{pmatrix}
	:=
	\begin{bmatrix} I_{n_J} & 0 & 0 \\ 0 & -b_N & -K \\ 0 & 1 & 0 \end{bmatrix} \begin{pmatrix} \hat{x} \\ \hat{u} \end{pmatrix}.
\end{equation*}
Then by \Cref{thm:augmentedsystem}, we have that \((\begin{psmallmatrix} x_J \\ x_N \end{psmallmatrix},u\bigr)\) is a solution of~\eqref{eq:dae} and
\begin{equation*}
	(E x)(0)
	= \begin{pmatrix} x_J(0) \\ N x_N(0) \end{pmatrix}
	= \begin{bmatrix} I_{n_J} & 0 \\ 0 & -K \end{bmatrix} \hat{x}(0)
	= F \hat{x}(0)
	= x^0.
\end{equation*}
So by \eqref{eq:costfunctionalequality}, we have
\begin{equation*}
	V_T(x^0) \leq J_T(x,u) = \hat{J}_T(\hat{x},\hat{u}) \leq \hat{V}_T(x^0)+\varepsilon.
\end{equation*}
As \(\varepsilon>0\) is arbitrary, we get
\begin{equation*}
	V_T(x^0) \leq \hat{V}_T(x^0).
\end{equation*}
To prove the reverse inequality, let \(x^0\in\consistencyspace\) and \(\varepsilon>0\) be arbitrary. Choose a solution~\((x,u)\) of~\eqref{eq:dae} such that \((E x)(0)=x^0\) and \(J_T(x,u)\leq V_T(x^0)+\varepsilon\). By \Cref{thm:augmentedsystem}, we have that
\begin{equation*}
	\hat{x} := \begin{pmatrix} x_J^\top, & u, & \dots, & u^{(\inputindex-1)} \end{pmatrix}^\top
\end{equation*}
solves the ordinary differential equation~\eqref{eq:xo} (or~\eqref{eq:augmentedsystem:one} for \(\inputindex=0\)).
Furthermore,
\begin{equation*}
	(N x_N)(0) = -N \begin{bmatrix} 0_{n_N\times n_J} & b_N & K \end{bmatrix} \begin{pmatrix} \hat{x}(0) \\ u^{(\inputindex)}(0) \end{pmatrix} = \begin{bmatrix} 0_{n_N\times n_J} & -K \end{bmatrix} \hat{x}(0),
\end{equation*}
therefore
\begin{equation*}
	F \hat{x}(0)
	= \begin{bmatrix} I_{n_J} & 0 \\ 0 & -K \end{bmatrix} \hat{x}(0)
	= \begin{pmatrix} x_J(0) \\ (N x_N)(0) \end{pmatrix} = (E x)(0).
\end{equation*}
Consequently
\begin{equation*}
	\hat{V}_T(x^0) \leq \hat{J}_T(\hat{x},u^{(\inputindex-1)}) = J_T(x,u) \leq V_T(x^0)+\varepsilon.
\end{equation*}
As for \(\varepsilon\to 0\), we obtain
\begin{equation*}
	\hat{V}_T(x^0) \leq V_T(x^0),
\end{equation*}
the assertion is proved.
\end{proof}

\begin{proof}[\textbf{Proof of~\Cref{thm:singularfeedback}}]
As~\eqref{eq:dae} is regular and therefore \((s E-A)\neq0_{\mathbb{R}[s]}\), the Sherman-Morrison-Woodbury formula (cf. Bernstein\cite{Bern09}, Fact~2.16.3)
\begin{equation*}
	\det((s E-A)+b k) = (1+k (s E-A)^{-1} b) \det(s E-A)
\end{equation*}
yields
\begin{equation}
	\det(s E-(A-b k))=0_{\mathbb{R}[s]} \ \Longleftrightarrow\ 1+k (s E-A)^{-1} b=0. \label{eq:singularsmw}
\end{equation}
Furthermore
\begin{align}
	k (s E-A)^{-1} b
	&= \begin{pmatrix} k_J & k_N \end{pmatrix} \begin{bmatrix} (s I_{n_J}-J)^{-1} & 0 \\ 0 & (s N-I_{n_N})^{-1} \end{bmatrix} \begin{pmatrix} b_J \\ b_N \end{pmatrix} \notag \\
	&= k_J (s I_{n_J}-J)^{-1} b_J + k_N (s N-I_{n_N})^{-1} b_N \in \mathbb{R}(s). \label{eq:tfqwf}
\end{align}
By~Bernstein\cite{Bern09}, Equation~(4.4.23), we have
\begin{equation}
	(s I_{n_J}-J)^{-1} = \frac{\sum_{i=0}^{n_J-1} J_i s^i}{\det(s I_{n_J}-J)} \ 
	\text{ for } J_0,\dots,J_{n_J-2}\in\mathbb{R}^{n_J\times n_J} \text{ and } J_{n_J-1}=I_{n_J} \label{eq:resolvent}
\end{equation}
and the nilpotency of~\(N\) gives
\begin{equation}
	(s N-I_{n_N})^{-1} = -\sum_{i=0}^{n_N-1} N^i s^i. \label{eq:sniinv}
\end{equation}

\enquote{\(\Rightarrow\)}:
Suppose that \(\det(s E-(A-b k))=0_{\mathbb{R}[s]}\). Then~\eqref{eq:tfqwf} applied to~\eqref{eq:singularsmw}, using the two identities~\eqref{eq:resolvent} and~\eqref{eq:sniinv} yields
\begin{equation}
	-\det(s I_{n_J}-J) =  \sum_{i=0}^{n_J-1} k_J J_i b_J s^i - \sum_{i=0}^{n_N-1} k_N N^i b_N s^i\cdot\det(sI_{n_J}-J).
	\label{eq:feedbacksingularpolynom}
\end{equation}
Since
\begin{equation*}
	\det(s I_{n_J}-J) = s^{n_J}+\alpha_{n_J-1} s^{n_J-1}+\dots+\alpha_0 s^0\in\mathbb{R}[s]
\end{equation*}
we conclude from~\eqref{eq:feedbacksingularpolynom} that
\begin{equation*}
	\sum_{i=0}^{n_N-1} k_N N^i b_N s^i = k_N b_N
\end{equation*}
and therefore
\begin{equation}
	-1 = -k_N b_N, \quad
	k_N K = 0_{1\times\inputindex}. \label{eq:knnbnsingular}
\end{equation}
Substituting~\eqref{eq:knnbnsingular} in~\eqref{eq:feedbacksingularpolynom} now yields
\begin{equation*}
	\sum_{i=0}^{n_J-1} k_J J_i b_J s^i = 0.
\end{equation*}
This applied to~\eqref{eq:resolvent} shows the first of the necessary conditions in~\eqref{eq:feedbacksingular}. The second necessary condition was already shown in~\eqref{eq:knnbnsingular}.

\enquote{\(\Leftarrow\)}:
Substituting the necessary conditions in~\eqref{eq:feedbacksingular} into~\eqref{eq:tfqwf} together with~\eqref{eq:sniinv} directly shows
\begin{equation*}
	k (sE-A)^{-1} b = -1
\end{equation*}
and \(\det(s E-(A-b k))=0_{\mathbb{R}[s]}\) follows from~\eqref{eq:singularsmw}. This concludes the proof.
\end{proof}

\begin{proof}[\textbf{Proof of~\Cref{thm:optimalfeedback}}]
Define \(p:=\hat{r}^{-1} (\hat{b}^\top P(\infty) +\hat{h}^\top)\) and note that
\begin{equation*}
	\hat{k}_\alpha \begin{bmatrix} I_{\hat{n}} \\ -p \end{bmatrix} = \Bigl(\alpha\begin{pmatrix} p & 1 \end{pmatrix}+\begin{pmatrix} 0_{1\times n_J} & 1 & 0_{1\times\inputindex} \end{pmatrix}\Bigr) \begin{bmatrix} I_{\hat{n}} \\ -p \end{bmatrix} = \begin{pmatrix} 0_{1\times n_J} & 1 & 0_{1\times\inputindex} \end{pmatrix}.
\end{equation*}
Let \((x^*,u^*)\) be an optimal trajectory for~\eqref{eq:daeocp} according to \Cref{thm:optimalcontrol}. Furthermore, let
\begin{equation*}
	\hat{x}^* = \begin{pmatrix} (x_J^*)^\top, & u^*, & \dots, & {u^*}^{(\inputindex-1)} \end{pmatrix}^\top, \quad
	\hat{u}^* = {u^*}^{(\inputindex)}
\end{equation*}
be the corresponding solution of~\eqref{eq:augmentedsystem} according to~\Cref{thm:augmentedsystem}. By \Cref{thm:optimalvalue} and \Cref{thm:optimalcontrol} it follows that
\begin{align}
	u^* = \hat{x}^*_{n_J+1} = \begin{pmatrix} 0_{1\times n_J} & 1 & 0_{1\times\inputindex} \end{pmatrix} \hat{x}^* = \hat{k}_\alpha \begin{bmatrix} I_{\hat{n}} \\ -p \end{bmatrix} \hat{x}^* = \hat{k}_\alpha \begin{pmatrix} \hat{x}^*  \\ -p \hat{x}^* \end{pmatrix} = \hat{k}_\alpha \begin{pmatrix} \hat{x}^*  \\ \hat{u}^* \end{pmatrix}. \label{eq:uinxhat}
\end{align}
By~\eqref{eq:solxn-c}, we have
\begin{equation*}
	x
	= \begin{pmatrix} x_J^* \\ x_N^* \end{pmatrix} \overset{\mathrm{ae}}{=} \begin{bmatrix} I_{n_J} & 0 \\ 0 & -[b_N,\dots,N^{\inputindex} b_N] \end{bmatrix} \begin{pmatrix} x_J^* \\ u^* \\ \vdots \\ {u^*}^{(\inputindex)} \end{pmatrix}
	= G \begin{pmatrix} \hat{x}^* \\ \hat{u}^* \end{pmatrix}.
\end{equation*}
As in~\Cref{thm:consistentiv}\ref{item:kalmanrank}, it follows that \([b_N,\dots,N^{\inputindex} b_N]\) has full column rank, so \(G\) is left invertible and
\begin{equation*}
	\begin{pmatrix} \hat{x}^* \\ \hat{u}^* \end{pmatrix} \overset{\mathrm{ae}}{=} G^\dagger x^*(\cdot).
\end{equation*}
Substituting this into~\eqref{eq:uinxhat} yields that the optimal trajectory of~\eqref{eq:daeocp} fulfils
\begin{equation*}
	\tfrac{\mathrm{d}}{\mathrm{d} t} (E x^*)(t) = A x^*(t)+b u^*(t) = A x^*(t)+b \hat{k}_\alpha G^\dagger x^*(t) = (A+b \hat{k}_\alpha G^\dagger) x^*(t),
\end{equation*}
so~\(x^*\) is a solution of~\eqref{eq:optimalclosedloop}.

To ensure that~\(x^*\) is the only solution of~\eqref{eq:optimalclosedloop}, it remains to show that the closed-loop system~\eqref{eq:optimalclosedloop} is regular, i.\,e.\@ \(\det(s E-(A+b \hat{k}_\alpha G^\dagger))\neq0_{\mathbb{R}[s]}\). Partition \(-\hat{k}_\alpha G^\dagger=:(k_J,k_N)\in\mathbb{R}^{1\times n_J}\times\mathbb{R}^{1\times n_N}\) and note that
\begin{align*}
	k_N K
	&= (k_J,k_N) \begin{pmatrix} 0_{n_J\times1} \\ b_N \end{pmatrix}
	= -\hat{k}_\alpha \, G^\dagger \begin{pmatrix} 0_{n_J\times\inputindex} \\ K \end{pmatrix}
	= \hat{k}_\alpha\, G^\dagger G \begin{pmatrix} 0_{(n_J+1)\times\inputindex} \\ I_{\inputindex} \end{pmatrix} \\
	&= \begin{pmatrix} p & 1 \end{pmatrix} \begin{pmatrix} 0_{(n_J+1)\times\inputindex} \\ I_{\inputindex} \end{pmatrix} \\
	&= \begin{pmatrix} p_{n_J+2} & \dots & p_{\hat{n}} & 1 \end{pmatrix}
	\neq 0,
\end{align*}
and regularity follows by~\Cref{thm:singularfeedback}.
\end{proof}

\end{document}